\documentclass[11pt]{article}%
\usepackage{amsmath}
\usepackage{amsfonts}
\usepackage{amssymb}
\usepackage{graphicx}
\usepackage{a4wide}%
\newtheorem{theorem}{Theorem}

\newtheorem{lemma}[theorem]{Lemma}

\newtheorem{proposition}[theorem]{Proposition}
\newtheorem{remark}[theorem]{Remark}

\begin{document}

\title{On D.Y. Gao and X. Lu paper \textquotedblleft On the extrema of a nonconvex
functional with double-well potential in 1D"}
\author{Constantin Z\u{a}linescu\thanks{Faculty of Mathematics, University Al.\ I.
Cuza, Bd.\ Carol I, Nr. 11, 700506 Ia\c{s}i, Romania, e-mail:
\texttt{zalinesc@uaic.ro}. }}
\date{}
\maketitle

\begin{abstract}
The aim of this paper is to discuss the main result in the paper by
D.Y.~Gao and X.~Lu [\emph{On the extrema of a nonconvex functional
with double-well potential in 1D}, Z. Angew. Math. Phys. (2016)
67:62]. More precisely we provide a detailed study of the problem
considered in that paper, pointing out the importance of the norm on
the space $C^{1}[a,b]$; because no norm (topology) is mentioned on
$C^{1}[a,b]$ we look at it as being a subspace of $W^{1,p}(a,b)$ for
$p\in [1,\infty]$ endowed with its usual norm. We show that the
objective function has not local extrema with the mentioned
constraints for $p\in [1,4)$, and has (up to an additive constant)
only a local maximizer for $p=\infty$, unlike the conclusion of the
main result of the discussed paper where it is mentioned that there
are (up to additive constants) two local minimizers and a local
maximizer. We also show that the same conclusions are valid for the
similar problem treated in the preprint by X.~Lu and D.Y.~Gao
[\emph{On the extrema of a nonconvex functional with double-well
potential in higher dimensions}, arXiv:1607.03995].

\end{abstract}

\section{Introduction}

For a faithful presentation of the the problem and result discussed
in \cite{GaoLu:16} we quote from this paper:

\medskip

\textquotedblleft The fourth-order polynomial defined by

$H(x):=\nu/2(1/2x^{2}-\lambda)^{2}$, where $x\in\mathbb{R},$
$\nu,\lambda$ are positive constants (1)

\noindent is the well-known Landau's second-order free energy, each of its
local minimizers represents a possible phase state of the material, while each
local maximizer characterizes the critical conditions that lead to the phase transitions.

...

The purpose of this paper is to find the extrema of the following nonconvex
total potential energy functional in 1D,

$I[u]:=\left(  \int_{a}^{b}H\left(  \frac{du}{dx}\right)  -fu\right)  dx$. (2)

\noindent The function $f\in C[a,b]$ satisfies the normalized balance condition

$\int_{a}^{b}f(x)dx=0$, (3)

\noindent and

there exists a unique zero root for $f$ in $[a,b]$. (4)

\noindent Moreover, its $L^{1}$-norm is sufficiently small such that

$\left\Vert f\right\Vert _{L^{1}(a,b)}<2\lambda\nu\sqrt{2\lambda}/(3\sqrt{3}%
)$. (5)

The above assumption is reasonable since large $\left\Vert f\right\Vert
_{L^{1}(a,b)}$ may possibly lead to instant fracture, which is represented by
nonsmooth solutions. The deformation $u$ is subject to the following two constraints,

$u\in C^{1}[a,b]$, (6)

\smallskip$\frac{du}{dx}(a)=\frac{du}{dx}(b)=0$. (7)

...

Before introducing the main result, we denote

$F(x):=-\int_{a}^{x}f(\rho)d\rho,~~x\in\lbrack a,b]$.

\noindent Next, we define a polynomial of third order as follows,

$E(y):=2y^{2}(\lambda+y/\nu),~~y\in\lbrack-\nu\lambda,+\infty)$.

\noindent Furthermore, for any $A\in\lbrack0,8\lambda^{3}\nu^{2}/27)$,

$E_{3}^{-1}(A)\leq E_{2}^{-1}(A)\leq E_{1}^{-1}(A)$

\noindent stand for the three real-valued roots for the equation $E(y)=A$.

At the moment, we would like to introduce the main theorem.

\textbf{Theorem 1.1.} For any function $f\in C[a,b]$ satisfying (3)--(5), one
can find the local extrema for the nonconvex functional (2).

\textbullet\ For any $x\in\lbrack a,b]$, $\overline{u}_{1}$ defined below is a
local minimizer for the nonconvex functional (2),

\noindent$\overline{u}_{1}(x)=\int_{a}^{x}F(\rho)/E_{1}^{-1}(F^{2}(\rho
))d\rho+C_{1},~~\forall C_{1}\in\mathbb{R}$. (9)

\textbullet\ For any $x\in\lbrack a,b]$, $\overline{u}_{2}$ defined below is a
local minimizer for the nonconvex functional (2),

$\overline{u}_{2}(x)=\int_{a}^{x}F(\rho)/E_{2}^{-1}(F^{2}(\rho))d\rho
+C_{2},~~\forall C_{2}\in\mathbb{R}$. (10)

\textbullet\ For any $x\in\lbrack a,b]$, $\overline{u}_{3}$ defined below is a
local maximizer for the nonconvex functional (2),

$\overline{u}_{3}(x)=\int_{a}^{x}F(\rho)/E_{3}^{-1}(F^{2}(\rho))d\rho
+C_{3},~~\forall C_{3}\in\mathbb{R}$. (11)\textquotedblright

\medskip

As mentioned in \cite{GaoLu:16}, in getting the above result the
authors use ``the canonical duality method".

\medskip

Let us observe from the beginning that nothing is said about the norm (and the
corresponding topology) on $C^{1}[a,b]$ when speaking about local extrema
(minimizers or maximizers).

In the following we discuss a slightly more general problem and
compare our conclusions with those of Theorem 1.1 in
\cite{GaoLu:16}. We don't analyze the method by which the
conclusions in Theorem 1.1 of \cite{GaoLu:16} are obtained even if
this is worth being done. Similar problems are considered by Gao and
Ogden in \cite{GaoOgd:08} and \cite{GaoOgd:08z} which are discussed
by Voisei and Z\u{a}linescu in \cite{VoiZal:11} and \cite{VoiZal:12},
respectively.

More precisely consider $\theta\in C[a,b]$ such that $\theta(x)>0$ for
$x\in\lbrack a,b]$, the polynomial $H$ defined by $H(y):=\tfrac{1}{2}%
(\tfrac{1}{2}y^{2}-\lambda)^{2}$ with $\lambda>0$, and the function
\[
J:=J_{f}:C^{1}[a,b]\rightarrow\mathbb{R},\quad J_{f}(u):=\int_{a}^{b}%
\theta\cdot\left(  H\circ u^{\prime}-fu\right)  ,
\]
where, $\int_{a}^{b}h$ denotes the Riemann integral
$\int_{a}^{b}h(x)dx$ of the function $h:[a,b]\rightarrow\mathbb{R}$
(when it exists). Of course, taking $\theta$ the constant function
$\nu$ $(>0)$ and replacing $f$ by $\nu^{-1}f$ we get the functional
$I$ considered in \cite{GaoLu:16}.

Let us set
\begin{align*}
X  &  :=C_{0}[a,b]:=\{v\in C[a,b]\mid v(a)=v(b)=0\},\\
Y  &  :=C_{1,0}[a,b]:=\{u\in C^{1}[a,b]\mid u^{\prime}:=du/dx\in C_{0}[a,b]\}.
\end{align*}
Of course $X$ is a linear subspace of $C[a,b]$; it is even a closed subspace
(and so a Banach space) if $C[a,b]$ is endowed with the supremum norm
$\left\Vert \cdot\right\Vert _{\infty}$. Clearly, other norms could be
considered on $X$.

Observe that the function $F$ defined in \cite{GaoLu:16} (and quoted
above) is in $C^{1}[a,b]\cap X$ with $F^{\prime}:=dF/dx=-f$.
Moreover, condition (5) implies that $\left\Vert F\right\Vert
_{\infty}<2\lambda\sqrt{2\lambda }/(3\sqrt{3})=(2\lambda/3)^{3/2}$
because
\[
\left\vert F(x)\right\vert =\left\vert \int_{a}^{x}f(\xi)d\xi\right\vert
\leq\int_{a}^{x}\left\vert f(\xi)\right\vert d\xi\leq\int_{a}^{b}\left\vert
f(\xi)\right\vert d\xi=\left\Vert f\right\Vert _{L^{1}(a,b)}.
\]
Furthermore, condition (4) implies that $F(x)>0$ for $x\in(a,b)$, or
$F(x)<0$ for $x\in(a,b)$.

For $u\in Y$ and $v:=u^{\prime}$ we have that
\begin{equation}
\int_{a}^{b}uf=-\int_{a}^{b}uF^{\prime}=-\left.  u(x)F(x)\right\vert _{a}%
^{b}+\int_{a}^{b}u^{\prime}F=\int_{a}^{b}vF. \label{r-gl1}%
\end{equation}
Using this fact, for $u$ satisfying the constraints (6) and (7), and
$v:=u^{\prime}$, one has
\[
J(u)=\int_{a}^{b}\theta\left(  H\circ v-Fv\right)  =:K(v).
\]

\section{Study of local extrema of the function $K$}

As mentioned above, in the sequel
$H:\mathbb{R}\rightarrow\mathbb{R}$ is defined by
$H(y):=\tfrac{1}{2}\left(  \tfrac{1}{2}y^{2}-\lambda\right)  ^{2}$
with $\lambda>0$, $\theta\in C[a,b]$ is such that
$\mu:=\min_{x\in\lbrack a,b]}\theta(x)>0$; moreover $F\in
C^{1}[a,b]\cap X$ is such that $F(x)\ne 0$ for $x\in (a,b)$ and
$\left\Vert F\right\Vert _{\infty}<(2\lambda/3)^{3/2}$.

Our first purpose is to find the local extrema of
\begin{equation}
K:=K_{F}:X\rightarrow\mathbb{R},\quad K_{F}(v):=\int_{a}^{b}\theta\cdot\left(
H\circ v-Fv\right)  \label{r-K}%
\end{equation}
on $X=C_{0}[a,b]$ endowed with the norm $\left\Vert \cdot\right\Vert _{p}$,
where $p\in\lbrack1,\infty]$.

First we study the Fr\'{e}chet and G\^{a}teaux differentiability of $K.$

\begin{lemma}
\label{l1}Let $g\in C[a,b]\setminus\{0\}$, $s\in\mathbb{N}^{\ast}%
\setminus\{1\}$ and $p\in\lbrack1,\infty]$. Then, with $h\in X,$
\[
\lim_{\left\Vert h\right\Vert _{p}\rightarrow0}\frac{1}{\left\Vert
h\right\Vert _{p}}\int_{a}^{b}gh^{s}=0\iff p\geq s.
\]

\end{lemma}

Proof. Set $\gamma:=\left\Vert g\right\Vert _{\infty}$ $(>0)$. For
$s<p<\infty$ and $h\in X$ we have that
\[
\left\vert \int_{a}^{b}gh^{s}\right\vert \leq\gamma\int_{a}^{b}\left\vert
h\right\vert ^{s}\cdot1\leq\gamma\left(  \int_{a}^{b}\left(  \left\vert
h\right\vert ^{s}\right)  ^{p/s}\right)  ^{s/p}\left(  \int_{a}^{b}%
1^{p/(p-s)}\right)  ^{(p-s)/p},
\]
and so
\begin{equation}
\left\vert \int_{a}^{b}gh^{s}\right\vert \leq\gamma(b-a)^{^{(p-s)/p}%
}\left\Vert h\right\Vert _{p}^{s}\quad\forall h\in X.\label{r-gl7}%
\end{equation}
The above inequality is true also, as easily seen, for $p=s$ and $p=\infty$
(setting $(p-s)/p:=1$ in the former case); from it we get $\lim_{\left\Vert
h\right\Vert _{p}\rightarrow0}\frac{1}{\left\Vert h\right\Vert _{p}}\int
_{a}^{b}gh^{s}=0$ because $s>1.$

Assume now that $p<s$. Since $g\in C[a,b]\setminus\{0\}$, there
exist $\delta>0$, and $a^{\prime},b^{\prime}\in\lbrack a,b]$ with
$a^{\prime
}<b^{\prime}$ such that $g(x)\geq\delta$ for $x\in\lbrack a^{\prime}%
,b^{\prime}]$ or $g(x)\leq-\delta$ for $x\in\lbrack a^{\prime},b^{\prime}]$.
Doing a translation, we suppose that $a^{\prime}=0$. For $n\in\mathbb{N}%
^{\ast}$ with $n\geq n_{0}$ $(\geq2/b^{\prime})$ consider
\begin{equation}
h_{n}(x):=\left\{
\begin{array}
[c]{ll}%
\alpha_{n}x & \text{if }x\in\lbrack0,1/n],\\
\alpha_{n}(2/n-x) & \text{if }x\in(1/n,2/n),\\
0 & \text{if }x\in\lbrack a,0)\cup\lbrack2/n,b],
\end{array}
\right.  \label{r-gl8}%
\end{equation}
with $\alpha_{n}:=n^{1+\gamma/p}>0$, where
$\frac{p-1}{s-1}<\gamma<1$.
Clearly, $h_{n}\in X=C_{0}[a,b]$. In this situation%
\[
\left\vert \int_{a}^{b}gh_{n}^{s}\right\vert =\int_{0}^{2/n}\left\vert
g\right\vert h_{n}^{s}\geq2\delta\int_{0}^{1/n}(\alpha_{n}x)^{s}%
dx=2\delta\alpha_{n}^{s}\frac{1}{s+1}\frac{1}{n^{s+1}}=\frac{2\delta}%
{s+1}n^{\frac{s\gamma-p}{p}},
\]
while a similar argument gives
\[
\left\Vert h_{n}\right\Vert _{p}=\left(  2\alpha_{n}^{p}\frac{1}{p+1}\frac
{1}{n^{p+1}}\right)  ^{1/p}=\left(  \frac{2}{p+1}\right)  ^{1/p}%
n^{\frac{\gamma-1}{p}}\rightarrow0.
\]
On the other hand,
\[
\frac{1}{\left\Vert h_{n}\right\Vert _{p}}\left\vert \int_{a}^{b}gh_{n}%
^{s}\right\vert \geq\frac{2\delta}{s+1}\left(  \frac{p+1}{2}\right)
^{1/p}n^{\frac{\gamma(s-1)-(p-1)}{p}}\rightarrow\infty,
\]
which proves our assertion. The proof is complete. \hfill$\square$

\begin{proposition}
\label{p1}Let $X=C_{0}[a,b]$ be endowed with the norm $\left\Vert
\cdot\right\Vert _{p}$, where $p\in\lbrack1,\infty]$. Then $K$ is
G\^{a}teaux differentiable; moreover, for $v\in X$, $K$ is Fr\'{e}chet
differentiable at $v$ if and only if $p\geq4$.
\end{proposition}

Proof. Let us set $g_{2}:=\tfrac{1}{2}\theta\left(  \tfrac{3}{2}v^{2}%
-\lambda\right)  $, $g_{3}:=\tfrac{1}{2}\theta v$ and $g_{4}:=\tfrac{1}%
{8}\theta$; of course, $g_{2},g_{3},g_{4}\in C[a,b]$. Set also $\beta
:=\max\{\left\Vert g_{2}\right\Vert _{\infty},\left\Vert g_{3}\right\Vert
_{\infty}\}$.

Observe that for all $v,h\in X$ we have that
\begin{equation}
K(v+h)=K(v)+\int_{a}^{b}\theta\left[  v(\tfrac{1}{2}v^{2}-\lambda)-F\right]
h+\int_{a}^{b}\tfrac{1}{2}\theta\left(  \tfrac{3}{2}v^{2}-\lambda\right)
h^{2}+\int_{a}^{b}\tfrac{1}{2}\theta vh^{3}+\int_{a}^{b}\tfrac{1}{8}\theta
h^{4}. \label{r-gl2}%
\end{equation}
For $v\in X$ consider
\begin{equation}
T_{v}:X\rightarrow\mathbb{R},\quad T_{v}(h):=\int_{a}^{b}\theta\left[
v(\tfrac{1}{2}v^{2}-\lambda)-F\right]  h\quad(h\in X). \label{r-gl6}%
\end{equation}
Clearly, $T_{v}$ is a linear operator; $T_{v}$ is also continuous for every
$p\in\lbrack1,\infty]$. Indeed, setting $\gamma_{v}:=\left\Vert \theta\left[
v(\tfrac{1}{2}v^{2}-\lambda)-F\right]  \right\Vert _{\infty}\in\mathbb{R}_{+}$
we have that
\[
\left\vert T_{v}(h)\right\vert \leq\gamma_{v}\int_{a}^{b}\left\vert
h\right\vert \leq\gamma_{v}\left\Vert h\right\Vert _{p}\cdot\left\Vert
1\right\Vert _{p^{\prime}}=\gamma_{v}(b-a)^{1/p^{\prime}}\left\Vert
h\right\Vert _{p}\quad\forall h\in X
\]
for $p,p^{\prime}\in\lbrack1,\infty]$ with $p^{\prime}$ the
conjugate of $p,$ that is $p^{\prime}:=p/(p-1)$ for
$p\in(1,\infty)$, $p^{\prime}:=\infty$ for $p=1$ and $p^{\prime}:=1$
for $p=\infty$. Hence $T_{v}$ is continuous.

Let $p\in\lbrack1,\infty]$ and $v\in X$ be fixed. Using (\ref{r-gl2}) we have
that
\[
\left\vert \frac{K(v+h)-K(v)-T_{v}(h)}{\left\Vert h\right\Vert _{p}%
}\right\vert \leq\frac{1}{\left\Vert h\right\Vert _{p}}\left(  \left\vert
\int_{a}^{b}g_{2}h^{2}\right\vert +\left\vert \int_{a}^{b}g_{3}h^{3}%
\right\vert +\left\vert \int_{a}^{b}g_{4}h^{4}\right\vert \right)
\]
for $h\neq0$. Using Lemma \ref{l1} for $p\geq4$, we obtain that $\lim
_{\left\Vert h\right\Vert _{p}\rightarrow0}\frac{K(v+h)-K(v)-T_{v}%
(h)}{\left\Vert h\right\Vert _{p}}=0$. Hence $K$ is Fr\'{e}chet differentiable
at $v.$

Assume now that $p<4$. Using again (\ref{r-gl2}) we have that%
\[
K(v+h)-K(v)-T_{v}(h)\geq\tfrac{\mu}{8}\int_{a}^{b}h^{4}-\beta\int_{a}%
^{b}\left\vert h\right\vert ^{3}-\beta\int_{a}^{b}h^{2}\quad\forall h\in X.
\]
Take $a=a^{\prime}=0<b^{\prime}=b$ (possible after a translation), $\alpha
_{n}:=n^{1+\gamma/p}$ with $\frac{3}{s-1}<\gamma<1$ and $h:=h_{n}$ defined by
(\ref{r-gl8}). Using the computations from the proof of Lemma \ref{l1}, we
get
\begin{equation}
\int_{a}^{b}\left\vert h_{n}\right\vert ^{s}=2\int_{0}^{1/n}(\alpha_{n}%
x)^{s}dx=\frac{2}{s+1}n^{\frac{s\gamma-p}{p}},\quad\left\Vert h_{n}\right\Vert
_{p}=\left(  \frac{2}{p+1}\right)  ^{1/p}\frac{1}{n^{(1-\gamma)/p}}%
\rightarrow0, \label{r-gl11}%
\end{equation}
whence
\begin{align*}
\frac{K(v+h_{n})-K(v)-T_{v}(h_{n})}{\left\Vert h_{n}\right\Vert _{p}}  &
\textstyle\geq\left(  \frac{p+1}{2}\right)  ^{1/p}n^{\frac{1-\gamma}{p}%
}\left(  \tfrac{\mu}{8}\cdot\tfrac{2}{5}n^{\frac{4\gamma-p}{p}}-\tfrac{2}%
{4}\beta_{3}n^{\frac{3\gamma-p}{p}}-\tfrac{2}{3}\beta_{2}n^{\frac{2\gamma
-p}{p}}\right) \\
&  \textstyle=\left(  \frac{p+1}{2}\right)  ^{1/p}n^{\frac{1-p+3\gamma}{p}%
}\left(  \tfrac{\mu}{20}-\tfrac{1}{2}\beta_{3}n^{-\frac{\gamma}{p}}-\tfrac
{2}{3}\beta_{2}n^{-\frac{2\gamma}{p}}\right)  \rightarrow\infty.
\end{align*}
This shows that $K\ $is not Fr\'{e}chet differentiable at $v.$

Because $K:(X,\left\Vert \cdot\right\Vert _{\infty})\rightarrow\mathbb{R}$ is
Fr\'{e}chet differentiable at $v\in X$, it follows that
\begin{equation}
\lim_{t\rightarrow0}\frac{K(v+th)-K(v)}{t}=T_{v}(h)\in\mathbb{R}\quad\forall
h\in X. \label{r-gl5}%
\end{equation}
Because $T_{v}:(X,\left\Vert \cdot\right\Vert
_{p})\rightarrow\mathbb{R}$ is linear and continuous, it follows
that $K$ is G\^{a}teaux differentiable at $v$ for every
$p\in\lbrack1,\infty]$ with $\nabla K(v)=T_{v}$. \hfill$\square$

\medskip

We consider now the problem of finding the stationary points of $K$,
that is those points $v\in X$ with $T_{v}=0.$

\begin{proposition}
\label{p2}The functional $K$ has only one stationary point
$\overline{v}$. More precisely, for each $x\in\lbrack a,b]$,
$\overline{v}(x)$ is the unique solution from
$(-\sqrt{2\lambda/3},\sqrt{2\lambda/3})$ of the equation
$z(\tfrac{1}{2}z^{2}-\lambda)=F(x)$.
\end{proposition}

Proof. Assume that $v\in X$ is stationary; hence
$T_{v}h=\int_{a}^{b}Vh=0$ for every $h\in X$, where $V:=\theta
v(\tfrac{1}{2}v^{2}-\lambda)-F$ $(\in X\subset C[a,b])$. We claim
that $V=0$. In the contrary case, since $V$ is continuous, there
exists $x_{0}\in(a,b)$ with $V(x_{0})\neq0$. Suppose that
$V(x_{0})>0$. By the continuity of $V$ there exist $a^{\prime},b^{\prime}%
\in\mathbb{R}$ such that $a<a^{\prime}<x_{0}<b^{\prime}<b$ and $V(x)>0$ for
every $x\in\lbrack a^{\prime},b^{\prime}]$. Take
\[
\overline{h}:[a,b]\rightarrow\mathbb{R},\quad\overline{h}(x):=\left\{
\begin{array}
[c]{ll}%
\frac{x-a^{\prime}}{b^{\prime}-a^{\prime}} & \text{if }x\in(a^{\prime}%
,\tfrac{1}{2}(a^{\prime}+b^{\prime})],\\
\frac{b^{\prime}-x}{b^{\prime}-a^{\prime}} & \text{if }x\in(\tfrac{1}%
{2}(a^{\prime}+b^{\prime}),b^{\prime}],\\
0 & \text{if }x\in\lbrack a,a^{\prime}]\cup(b^{\prime},b].
\end{array}
\right.
\]
Then $\overline{h}\in X$ and $\overline{h}(x)>0$ for $x\in(a^{\prime
},b^{\prime})$. Since $0=\int_{a}^{b}V\overline{h}=\int_{a^{\prime}%
}^{b^{\prime}}V\overline{h}$ and $V\overline{h}$ is continuous and nonnegative
on $[a^{\prime},b^{\prime}]$ we obtain that $V(x)\overline{h}(x)=0$ for
$x\in\lbrack a^{\prime},b^{\prime}]$, and so $0=V(x_{0})\overline{h}(x_{0}%
)>0$. This contradiction shows that $V=0$. The proof in the case $V(x_{0})<0$
reduces to the preceding one replacing $V$ by $-V$. Hence
\begin{equation}
\theta v(\tfrac{1}{2}v^{2}-\lambda)=F~~\text{on~~}[a,b]. \label{r-gl3}%
\end{equation}

Consider the polynomial function $G:\mathbb{R}\rightarrow\mathbb{R}$ defined
by $G(z):=z\left(  \tfrac{1}{2}z^{2}-\lambda\right)  $. Then $G^{\prime
}(z)=\tfrac{3}{2}z^{2}-\lambda$ having the zeros $\pm\kappa$, where
\begin{equation}
\kappa:=\sqrt{2\lambda/3}. \label{r-k}%
\end{equation}
The behavior of $G$ is given in the table below.

%\smallskip

\begin{center}
{\footnotesize
\begin{tabular}
[c]{c|ccccccccccccccccc}%
$z$ & $-\infty$ &  & $-2\kappa$ &  & $-\sqrt{3}\kappa$ &  & $-\kappa$ &  & $0$
&  & $\kappa$ &  & $\sqrt{3}\kappa$ &  & $2\kappa$ &  & $\infty$\\\hline
$G^{\prime}(z)$ &  & $+$ &  & $+$ &  & $+$ & $0$ & $-$ &  & $-$ & $0$ & $+$ &
& $+$ &  & $+$ & \\\hline
$G(z)$ & $-\infty$ & $\nearrow$ & $-\kappa^{3}$ & $\nearrow$ & $0$ &
$\nearrow$ & $\kappa^{3}$ & $\searrow$ & $0$ & $\searrow$ & $-\kappa^{3}$ &
$\nearrow$ & $0$ & $\nearrow$ & $\kappa^{3}$ & $\nearrow$ & $\infty$%
\end{tabular}
}
\end{center}

This table shows that the equation $G(z)=A$ with $A\in(-\sqrt{8\lambda^{3}%
/27},\sqrt{8\lambda^{3}/27})$ has three real solutions, more precisely,
\begin{equation}
z_{1}(A)\in(-2\kappa,-\kappa),\quad z_{2}(A)\in(-\kappa,\kappa),\quad
z_{3}(A)\in(\kappa,2\kappa). \label{r-gl4}%
\end{equation}
Moreover, the mappings $z_{i}:(-\kappa^{3},\kappa^{3})\rightarrow\mathbb{R}$
are continuous with $z_{1}(0)=-\sqrt{3}\kappa$, $z_{2}(0)=0$, $z_{3}%
(0)=\sqrt{3}\kappa$. This shows that $z_{i}\circ F\in X$ if and only $i=2,$
and so the only solution in $X$ of the equation $v(\tfrac{1}{2}v^{2}%
-\lambda)=F$ is $\overline{v}:=z_{2}\circ F$. \hfill$\square$

\medskip

Let us analyze if $\overline{v}:=z_{2}\circ F$ is a local extremum of $K.$

\begin{proposition}
\label{p3}Let $\overline{v}\in X$ be the stationary point of $K$.
Then $\overline{v}:=z_{2}\circ F$ [with $z_{2}$ defined in
(\ref{r-gl4})] is a local maximizer for $K$ with respect to
$\left\Vert \cdot\right\Vert _{\infty }$, and $\overline{v}$ is not
a local extremum point of $K$ with respect to $\left\Vert
\cdot\right\Vert _{p}$ for $p\in\lbrack1,4).$
\end{proposition}

Proof. Let us consider first the case $p=\infty$. From (\ref{r-gl2}) we get
\begin{equation}
K(\overline{v}+h)-K(\overline{v})=\int_{a}^{b}\theta\left[  \tfrac{1}%
{2}\left(  \tfrac{3}{2}\overline{v}^{2}-\lambda\right)  +\tfrac{1}{2}%
\overline{v}h+\tfrac{1}{8}h^{2}\right]  h^{2}\quad\forall h\in X.
\label{r-gl10}%
\end{equation}
Since $F\in C[a,b]$, there exists some $x_{0}\in\lbrack a,b]$ such
that $\left\Vert F\right\Vert _{\infty}=\left\vert
F(x_{0})\right\vert <(2\lambda/3)^{3/2}$, and so $\left\vert
\overline{v}(x)\right\vert \leq\left\vert
\overline{v}(x_{0})\right\vert =:\gamma<\sqrt{2\lambda/3}$ for
$x\in\lbrack a,b]$. It follows that $\tfrac{1}{2}\left(  \tfrac{3}{2}%
\overline{v}^{2}-\lambda\right)  \leq\tfrac{1}{2}\left(  \tfrac{3}{2}%
\gamma^{2}-\lambda\right)  =:-\eta<\tfrac{1}{2}\left(  \tfrac{3}{2}%
\frac{2\lambda}{3}-\lambda\right)  =0$. Hence
\begin{equation}
\tfrac{1}{2}\left(  \tfrac{3}{2}\overline{v}^{2}-\lambda\right)  +\tfrac{1}%
{2}\overline{v}h+\tfrac{1}{8}h^{2}\leq-\eta+\tfrac{1}{2}\gamma\left\Vert
h\right\Vert _{\infty}+\tfrac{1}{8}\left\Vert h\right\Vert _{\infty}%
^{2}<0\quad\forall h\in X,~\left\Vert h\right\Vert _{\infty}<\varepsilon,
\label{r-gl9}%
\end{equation}
where $\varepsilon:=2\big(\sqrt{\gamma^{2}+2\eta}-\gamma\big)$. It follows
that $\overline{v}$ is a (strict) local maximizer of $K.$

Assume now that $p\in\lbrack1,4)$. Of course, there exists a sequence
$(h_{n})_{n\geq1}\subset X\setminus\{0\}$ such that $\left\Vert h_{n}%
\right\Vert _{\infty}\rightarrow0$. Taking into account
(\ref{r-gl9}), we have that $K(\overline{v}+h_{n})<K(\overline{v})$
for large $n$. Since $\left\Vert h_{n}\right\Vert _{p}\rightarrow0$,
$\overline{v}$ is not a local minimizer of $K$ with respect to
$\left\Vert \cdot\right\Vert _{p}$. In the proof of Proposition
\ref{p1} we found a sequence $(h_{n})_{n\geq1}\subset
X\setminus\{0\}$ such that $\left\Vert h_{n}\right\Vert
_{p}\rightarrow0$ and
$\left\Vert h_{n}\right\Vert _{p}^{-1}\left(  K(\overline{v}+h_{n}%
)-K(\overline{v})-T_{\overline{v}}h_{n}\right)  \rightarrow\infty$.
Since $T_{\overline{v}}=0$, we obtain that
$K(\overline{v}+h_{n})-K(\overline{v})>0$ for large $n$, proving
that $\overline{v}$ is not a local maximizer of $K$. Hence
$\overline{v}$ is not a local extremum point of $K$. \hfill$\square$

\medskip

We don't know if $\overline{v}$ is a local maximizer of $K$ for $p\in
\lbrack4,\infty)$; having in view (\ref{r-gl9}), surely, $\overline{v}$ is not
a local minimizer of $K.$

Proposition \ref{p3} shows the importance of the norm (and more generally, of
the topology) on a space when speaking about local extrema.

\medskip

Let us establish now the relations between the local extrema of $J$
with the constraints (6) and (7) in \cite{GaoLu:16}, that is local
extrema of $J$ restricted to $C_{1,0}[a,b]$, and the local extrema
of $K$ in the case in which $C^{1}[a,b]$ is endowed with the (usual)
norm defined by
\begin{equation}
\left\Vert u\right\Vert :=\left\Vert u\right\Vert _{\infty}+\left\Vert
u^{\prime}\right\Vert _{\infty}\quad(u\in C^{1}[a,b]), \label{r-gl13}%
\end{equation}
and $C_{0}[a,b]$ is endowed with the norm $\left\Vert \cdot\right\Vert
_{\infty}.$

\begin{proposition}
\label{p4}Consider the norm $\left\Vert \cdot\right\Vert $ (defined
in (\ref{r-gl13})) on $C^{1}[a,b]$ and the norm $\left\Vert
\cdot\right\Vert _{\infty}$ on $C_{0}[a,b]$. If $\overline{u}$ is a
local minimizer (maximizer) of $J$ on $C_{1,0}[a,b]$, then
$\overline{u}^{\prime}$ is a local minimizer (maximizer) of $K$.
Conversely, if $\overline{v}$ is a local minimizer (maximizer) of
$K$, then $\overline{u}\in C^{1}[a,b]$ defined by $\overline
{u}(x):=u_{0}+\int_{a}^{x}\overline{v}(\xi)d\xi$ for $x\in\lbrack
a,b]$ and a fixed $u_{0}\in\mathbb{R}$ is a local minimizer
(maximizer) of $J$ on $C_{1,0}[a,b].$
\end{proposition}

Proof. Assume that $\overline{u}$ is a local minimizer of $J$ on
$C_{1,0}[a,b]$; hence $\overline{u}\in C_{1,0}[a,b]$. It follows that there
exists $r>0$ such that $J(\overline{u})\leq J(u)$ for every $u\in
C_{1,0}[a,b]$ with $\left\Vert u-\overline{u}\right\Vert <r$. Set
$\overline{v}:=\overline{u}^{\prime}$ and take $v\in X=C_{0}[a,b]$ with
$\left\Vert v-\overline{v}\right\Vert _{\infty}<r^{\prime}:=r/(1+b-a)$. Define
$u:[a,b]\rightarrow\mathbb{R}$ by $u(x):=\overline{u}(a)+\int_{a}^{x}%
v(\xi)d\xi$ for $x\in\lbrack a,b]$. Then $u\in C_{1,0}[a,b]$ and
$u^{\prime }=v$. Since
$\overline{u}(x)=\overline{u}(a)+\int_{a}^{x}\overline{v}(\xi
)d\xi$, we get
\[
\left\Vert u-\overline{u}\right\Vert =\left\Vert u-\overline{u}\right\Vert
_{\infty}+\left\Vert u^{\prime}-\overline{u}^{\prime}\right\Vert _{\infty}%
\leq(b-a)\left\Vert v-\overline{v}\right\Vert _{\infty}+\left\Vert
v-\overline{v}\right\Vert _{\infty}<r^{\prime}(1+b-a)=r.
\]
Hence $K(\overline{v})=J(\overline{u})\leq J(u)=K(v)$. This shows that
$\overline{v}$ is a local minimizer for $K$.

Conversely, assume that $\overline{v}$ is a local minimizer for $K$.
Then there exists $r>0$ such that $K(\overline{v})\leq K(v)$ for
$v\in C_{0}[a,b]$ with $\left\Vert v-\overline{v}\right\Vert <r$,
and take $u_{0}\in\mathbb{R}$ and
$\overline{u}:[a,b]\rightarrow\mathbb{R}$ defined by $\overline
{u}(x):=u_{0}+\int_{a}^{x}\overline{v}(\xi)d\xi$ for $x\in\lbrack
a,b]$. Then $\overline{u}\in C_{1,0}[a,b]$. Consider $u\in
C_{1,0}[a,b]$ with $\left\Vert u-\overline{u}\right\Vert <r$, that
is
\[
\left\Vert u-\overline{u}\right\Vert _{\infty}+\left\Vert u^{\prime}%
-\overline{u}^{\prime}\right\Vert _{\infty}=\left\Vert u-\overline
{u}\right\Vert _{\infty}+\left\Vert u^{\prime}-\overline{v}\right\Vert
_{\infty}<r;
\]
then $\left\Vert u^{\prime}-\overline{v}\right\Vert _{\infty}<r$.
Since $u^{\prime}\in C_{0}[a,b]$, it follows that
$J(u)=K(u^{\prime})\geq K(\overline{v})=J(\overline{u})$, and so
$\overline{u}$ is a local minimizer of $J$ on $C_{1,0}[a,b]$. The
case of local maximizers for $J$ and $K$ is treated similarly.
\hfill$\square$

\medskip

Putting together Propositions \ref{p2}, \ref{p3} and \ref{p4} we get the next result.

\begin{theorem}
\label{t-gl}Consider the norm $\left\Vert \cdot\right\Vert $ (defined in
(\ref{r-gl13})) on $C^{1}[a,b]$ and the norm $\left\Vert \cdot\right\Vert
_{\infty}$ on $C_{0}[a,b]$. Let $\overline{u}\in C_{1,0}[a,b]$ and set
$\overline{v}:=\overline{u}^{\prime}$. Then the following assertions are equivalent:

\emph{(i)} $\overline{u}$ is a local maximum point of $J$ restricted to
$C_{1,0}[a,b].$

\emph{(ii)} $\overline{u}$ is a local extremum point of $J$ restricted to
$C_{1,0}[a,b].$

\emph{(iii)} $\overline{v}$ is a stationary point of $K.$

\emph{(iv)} $\overline{v}$ is a local extremum point of $K.$

\emph{(v)} $\overline{v}$ is a local maximum point of $K.$

\emph{(vi)} $\overline{v}=z_{2}\circ F$, where $z_{2}(A)$ is the
unique solution of the equation $z\left(
\tfrac{1}{2}z^{2}-\lambda\right)  =A$ in the interval
$(-\sqrt{2\lambda/3},\sqrt{2\lambda/3}]$ for $A\in(-(2\lambda
/3)^{3/2},(2\lambda/3)^{3/2}).$

\emph{(vii)} there exists $u_{0}\in\mathbb{R}$ such that $\overline
{u}(x)=u_{0}+\int_{a}^{x}z_{2}(F(\rho))d\rho$ for every $x\in\lbrack a,b].$
\end{theorem}

\section{Discussion of Theorem 1.1 from Gao and Lu's paper \cite{GaoLu:16}}

First of all, we think that in the formulation of \cite[Th.
1.1]{GaoLu:16}, \textquotedblleft local extrema for the nonconvex
functional (2)\textquotedblright\ must be replaced by
\textquotedblleft local extrema for the nonconvex functional (2)
with the constraints (6) and (7)\textquotedblright,
\textquotedblleft local minimizer for the nonconvex functional
(2)\textquotedblright\ must be replaced by \textquotedblleft local
minimizer for the nonconvex functional (2) with the constraints (6)
and (7)\textquotedblright\ (2 times), and \textquotedblleft local
maximizer for the nonconvex functional (2)\textquotedblright\ must
be replaced by \textquotedblleft local maximizer for the nonconvex
functional (2) with the constraints (6) and (7)\textquotedblright.
Below, we interpret \cite[Th. 1.1]{GaoLu:16} with these
modifications.

As pointed in Introduction, no norms are considered on the spaces
mentioned in \cite{GaoLu:16}. For this reason in Theorem \ref{t-gl}
we considered the usual norms on $C^{1}[a,b]$ and $C_{0}[a,b]$;
these norms are used in this discussion. Moreover, let
$\theta(x):=1$ for $x\in\lbrack a,b]$ in Theorem \ref{t-gl} and
$\nu=1$ in \cite[Th. 1.1]{GaoLu:16}. In the conditions of \cite[Th.
1.1]{GaoLu:16} $F(x)>0$ for $x\in(a,b)$ or $F(x)<0$ for $x\in(a,b).$
For the present discussion we take the case $F>0$ on $(a,b).$

Assume that the mappings
\begin{equation}
\rho\mapsto F(\rho)/E_{j}^{-1}\left(  F^{2}(\rho\right)  )=:v_{j}(\rho)
\label{r-vj}%
\end{equation}
[where \textquotedblleft$E_{3}^{-1}(A)\leq E_{2}^{-1}(A)\leq E_{1}^{-1}(A)$
stand for the three real-valued roots for the equation $E(y)=A$%
\textquotedblright\ with $E(y)=2y^{2}(y+\lambda)$ and $A\in\lbrack
0,8\lambda^{3}/27)$] are well defined for $\rho\in\{a,b\}$ [there are no
problems for $\rho\in(a,b)$].

If \cite[Th. 1.1]{GaoLu:16} is true, then $v_{1},v_{2},v_{3}\in
C_{0}[a,b]$; moreover, $v_{1}$ and $v_{2}$ are local minimizers of
$K$, and $v_{3}$ is a local maximizer of $K$. This is of course
\textbf{false} taking into account Theorem \ref{t-gl} because $K$
has not local minimizers.

Because $z_{2}\circ F$ is the unique local maximizer of $K$, we must have that
$v_{3}=z_{2}\circ F$. Let us see if this is true. Because $z_{i}(A)$ are
solutions of the equation $G(z)=A$ and $E_{j}^{-1}(A)$ are solutions of the
equation $E(y)=A$, we must study the relationships among these numbers.

First, the behavior of $E$ is given in the next table.

\begin{center}
{\small
\begin{tabular}
[c]{c|ccccccccccc}%
$y$ & $-\infty$ &  & $-\lambda$ &  & $-\frac{2\lambda}{3}$ &  & $0$ &  &
$\frac{\lambda}{3}$ &  & $\infty$\\\hline
$E^{\prime}(y)$ &  & $+$ &  & $+$ & $0$ & $-$ & $0$ & $+$ &  & $+$ & \\\hline
$E(y)$ & $-\infty$ & $\nearrow$ & $0$ & $\nearrow$ & $\frac{8\lambda^{3}}{27}$
& $\searrow$ & $0$ & $\nearrow$ & $\frac{8\lambda^{3}}{27}$ & $\nearrow$ &
$\infty$%
\end{tabular}
}
\end{center}

Secondly, for $y,z,A\in\mathbb{C}\setminus\{0\}$ such that $yz=A$ we have
that
\begin{equation}
G(z)=A\Leftrightarrow\frac{A}{y}\left(  \frac{1}{2}\frac{A^{2}}{y^{2}}%
-\lambda\right)  =A\Leftrightarrow2y^{2}(y+\lambda)=A^{2}\Leftrightarrow
E(y)=A^{2}. \label{r-gl12}%
\end{equation}

Analyzing the behavior of $G$ and $E$ (recall that $\kappa=\sqrt{2\lambda/3}%
$), and the relation $yz=A$ for $A\neq0$ (mentioned above), the correspondence
among the solutions of the equations $G(z)=A$ and $E(y)=A^{2}$ for
$A\in(0,(2\lambda/3)^{3/2})$ is:
\begin{equation}
z_{1}(A)=A/E_{2}^{-1}(A^{2}),\quad z_{2}(A)=A/E_{3}^{-1}(A^{2}),\quad
z_{3}(A)=A/E_{1}^{-1}(A^{2}) \label{r-ze}%
\end{equation}
for all $A\in(0,\left(  2\lambda/3\right)  ^{3/2})$. This shows that
only the third assertion of \cite[Th. 1.1]{GaoLu:16} is true (of
course, considering the norm defined in (\ref{r-gl13}) on
$C^{1}[a,b]$).

\section{Discussion of Theorem 1.1 from Lu and Gao's paper \cite{LuGao:16}}

A similar problem to that in \cite{GaoLu:16}, discussed above, is
considered in \cite{LuGao:16}. In the abstract of this paper one
finds:

\textquotedblleft In comparison with the 1D case discussed by D. Gao and R.
Ogden, there exists huge difference in higher dimensions, which will be
explained in the theorem".

More precisely, in \cite{LuGao:16} it is said:

\textquotedblleft In this paper, we consider the fourth-order polynomial
defined by

$H(|\vec{\gamma}|):=\nu/2\left(
1/2|\vec{\gamma}|^{2}-\lambda\right)  ^{2},$
$\vec{\gamma}\in\mathbb{R}^{n}$, $\nu,\lambda>0$ are constants,
$|\vec{\gamma }|^{2}=\vec{\gamma}\cdot\vec{\gamma}$.

...

The purpose of this paper is to find the extrema of the following nonconvex
total potential energy functional in higher dimensions,

(1) $I[u]:=\int_{\Omega}\left(  H(|\nabla u|)-fu\right)  dx,$

\noindent where $\Omega=$Int$\left\{  \mathbb{B}(O,R_{1})\setminus
\mathbb{B}(O,R_{2})\right\}  $, $R_{1}>R_{2}>0,$
$\mathbb{B}(O,R_{1})$ and $\mathbb{B}(O,R_{2})$ denote two open
balls with center $O$ and radii $R_{1}$ and $R_{2}$ in the Euclidean
space $\mathbb{R}^{n}$, respectively. \textquotedblleft
Int\textquotedblright\ denotes the interior points. In
addition, let $\Sigma_{1}:=\{x:|x|=R_{1}\}$, and $\Sigma_{2}:=\{x:|x|=R_{2}%
\}$, then the boundary $\partial\Omega=\Sigma_{1}\cup\Sigma_{2}$. The radially
symmetric function $f\in C(\overline{\Omega})$ satisfies the normalized
balance condition

(2) $\int_{\Omega}f(|x|)dx=0$,

\noindent and

(3) $f(|x|)=0$ if and only if $|x|=R_{3}\in(R_{2},R_{1})$.

\noindent Moreover, its $L^{1}$-norm is sufficiently small such that

(4) $\left\Vert f\right\Vert _{L^{1}(\Omega)}<4\lambda\nu R_{2}^{n-1}%
\sqrt{2\lambda\pi^{n}}/(3\sqrt{3}\Gamma(n/2)),$

\noindent where $\Gamma$ stands for the Gamma function. This assumption is
reasonable since large $\left\Vert f\right\Vert _{L^{1}(\Omega)}$ may possibly
lead to instant fracture. The deformation $u$ is subject to the following
three constraints,

(5) $u$ is radially symmetric on $\overline{\Omega}$,

(6) $u\in W^{1,\infty}(\Omega)\cap C(\overline{\Omega})$,

(7) $\nabla u\cdot\vec{n}=0$ on both $\Sigma_{1}$ and $\Sigma_{2}$,

\noindent where $\vec{n}$ denotes the unit outward normal on $\partial\Omega$.

By variational calculus, one derives a correspondingly nonlinear
Euler--Lagrange equation for the primal nonconvex functional, namely,

(8) $\operatorname*{div}\left(  \nabla H(|\nabla u|)\right)  +f=0$ in $\Omega$,

\noindent equipped with the Neumann boundary condition (7). Clearly, (8) is a
highly nonlinear partial differential equation which is difficult to solve by
the direct approach or numerical method [2, 15]. However, by the canonical
duality method, one is able to demonstrate the existence of solutions for this
type of equations.

...

Before introducing the main result, we denote

$F(r):=-1/r^{n}\int_{R_{2}}^{r}f(\rho)\rho^{n-1}d\rho.~~r\in\lbrack
R_{2},R_{1}]$.

\noindent Next, we define a polynomial of third order as follows,

$E(y):=2y^{2}(\lambda+y/\nu),~~y\in\lbrack-\nu\lambda,+\infty)$.

\noindent Furthermore, for any $A\in\lbrack0,8\lambda^{3}\nu^{2}/27)$,

$E_{3}^{-1}(A)\leq E_{2}^{-1}(A)\leq E_{1}^{-1}(A)$

\noindent stand for the three real-valued roots for the equation $E(y)=A$.

At the moment, we would like to introduce the theorem of multiple extrema for
the nonconvex functional (2).

\textbf{Theorem 1.1.} For any radially symmetric function $f\in C(\overline
{\Omega})$ satisfying (2)--(4), we have three solutions for the nonlinear
Euler--Lagrange equation (8) equipped with the Neumann boundary condition, namely

\textbullet\ For any $r\in\lbrack R_{2},R_{1}]$, $\overline{u}_{1}$ defined
below is a local minimizer for the nonconvex functional (2),

\noindent$(9)\quad\overline{u}_{1}(\left\vert x\right\vert )=\overline{u}%
_{1}(r):=\int_{R_{2}}^{r}F(\rho)\rho/E_{1}^{-1}(F^{2}(\rho)\rho^{2}%
)d\rho+C_{1},~~\forall C_{1}\in\mathbb{R}$.

\textbullet\ For any $r\in\lbrack R_{2},R_{1}]$, $\overline{u}_{2}$ defined
below is a local minimizer for the nonconvex functional (2) in 1D. While for
the higher dimensions $n\geq2$, $\overline{u}_{2}$ is not necessarily a local
minimizer for (2) in comparison with the 1D case.

\noindent$(10)\quad\overline{u}_{2}(\left\vert x\right\vert )=\overline{u}%
_{2}(r):=\int_{R_{2}}^{r}F(\rho)\rho/E_{2}^{-1}(F^{2}(\rho)\rho^{2}%
)d\rho+C_{2},~~\forall C_{2}\in\mathbb{R}$.

\textbullet\ For any $r\in\lbrack R_{2},R_{1}]$, $\overline{u}_{3}$ defined
below is a local maximizer for the nonconvex functional (2),

\noindent$(11)\quad\overline{u}_{3}(\left\vert x\right\vert )=\overline{u}%
_{3}(r):=\int_{R_{2}}^{r}F(\rho)\rho/E_{3}^{-1}(F^{2}(\rho)\rho^{2}%
)d\rho+C_{3},~~\forall C_{3}\in\mathbb{R}$.

...

In the final analysis, we apply the canonical duality theory to prove Theorem
1.1.\textquotedblright

\medskip

First, observe that one must have (1) instead of (2) just before the
statement of \cite[Th.~1.1]{LuGao:16}, as well as in its statement,
excepting for (2)--(4). Secondly, (even from the quoted texts) one
must observe that the wording in \cite{GaoLu:16} and \cite{LuGao:16}
is almost the same; the mathematical part is very, very similar,
too.

\medskip

To avoid any confusion, in the sequel the Euclidian norm on $\mathbb{R}^{n}$
will be denoted by $\left\vert \cdot\right\vert _{n}$ instead of $\left\vert
\cdot\right\vert .$

Remark that it is said $f\in C(\overline{\Omega})$, which implies
$f$ is applied to elements $x\in\overline{\Omega}$, while a line
below one considers $f(\left\vert x\right\vert )$ (that is
$f(\left\vert x\right\vert _{n})$ with our notation); because the
(Euclidean) norm $\left\vert x\right\vert _{n}$ of
$x\in\overline{\Omega}$ belongs to $[R_{2},R_{1}]$, writing
$f(\left\vert x\right\vert )$ shows that
$f:[R_{2},R_{1}]\rightarrow\mathbb{R}$. Of course, these create
ambiguities. Probably the authors wished to say that a function
$g:\overline{\Omega}\rightarrow\mathbb{R}$ is radially symmetric if
there exists $\psi:[R_{2},R_{1}]\rightarrow\mathbb{R}$ such that
$g(x)=\psi (\left\vert x\right\vert _{n})$ for every
$x\in\overline{\Omega}$, that is $g=\psi\circ\left\vert
\cdot\right\vert _{n}$ on $\overline{\Omega}$; observe that $\psi$
is continuous if and only if $\psi\circ\left\vert \cdot\right\vert
_{n}$ is continuous. Because also the functions $u$ in the
definition of $I$
are asked to be radially symmetric on $\overline{\Omega}$ (see \cite[(5)]%
{LuGao:16}), it is useful to observe that for a Riemann integrable
function $\psi:[R_{2},R_{1}]\rightarrow\mathbb{R}$, using the usual
spherical change of variables, we have that
\begin{equation}
\int_{\Omega}\left(  \psi\circ\left\vert \cdot\right\vert _{n}\right)
(x)dx=\frac{2\pi^{n/2}}{\Gamma(n/2)}\cdot\int_{R_{2}}^{R_{1}}r^{n-1}%
\psi(r)dr=\gamma_{n}\int_{R_{2}}^{R_{1}}\theta\psi, \label{r-int}%
\end{equation}
where
\begin{equation}
\gamma_{n}:=\frac{2\pi^{n/2}}{\Gamma(n/2)},\text{~~and~~}\theta:[R_{2}%
,R_{1}]\rightarrow\mathbb{R},~~\theta(r):=r^{n-1}. \label{r-gn-teta}%
\end{equation}

So, in the sequel we consider that $f:[R_{2},R_{1}]\rightarrow\mathbb{R}$ is
continuous. Condition \cite[(2)]{LuGao:16} becomes $\int_{R_{2}}^{R_{1}}\theta
f=0$ [for the definition of $\theta$ see (\ref{r-gn-teta})], condition
\cite[(3)]{LuGao:16} is equivalent to the existence of a unique $R_{3}%
\in(R_{2},R_{1})$ such that $f(R_{3})=0$ (that is $(\theta f)(R_{3})=0$),
while condition \cite[(4)]{LuGao:16} is equivalent to $\big\Vert\theta
f\big\Vert_{L^{1}[R_{2},R_{1}]}<\nu R_{2}^{n-1}(2\lambda/3)^{3/2}.$

Moreover, condition \cite[(5)]{LuGao:16} is equivalent to the
existence of $\upsilon:[R_{2},R_{1}]\rightarrow\mathbb{R}$ such that
$u=\upsilon \circ\left\vert \cdot\right\vert _{n}$, while the
condition $u\in C(\overline{\Omega})$ is equivalent to $\upsilon\in
C[R_{2},R_{1}].$

Which is the meaning of $\nabla u(x)$ in condition \cite[(7)]{LuGao:16} for
$u\in W^{1,\infty}(\Omega)$ and $x\in\Sigma_{1}$ (or $x\in\Sigma_{2}$)? For
example, let us consider $\upsilon:[1,3]\rightarrow\mathbb{R}$ defined by
$\upsilon(t):=(t-1)^{2}\sin\frac{1}{t-1}$ for $t\in(1,2)$. Is $u:=\upsilon
\circ\left\vert \cdot\right\vert _{n}$ in $W^{1,\infty}(\Omega)$ for
$R_{2}:=1$ and $R_{1}:=2?$ If YES, which is $\nabla u(x)$ for $x\in
\mathbb{R}^{n}$ with $\left\vert x\right\vert _{n}=1?$

Let us assume that $\upsilon\in
C^{1}(R_{2}-\varepsilon,R_{1}+\varepsilon)$ for some
$\varepsilon\in(0,R_{2})$ and take $u:=\upsilon\circ\left\vert
\cdot\right\vert _{n}$. Then clearly $u\in C^{1}(\Delta)$, where
$\Delta:=\{x\in\mathbb{R}^{n}\mid\left\vert x\right\vert _{n}\in
(R_{2}-\varepsilon,R_{1}+\varepsilon)\}$, and
\begin{equation}
\nabla u(x)=\upsilon^{\prime}(\left\vert x\right\vert _{n})\left\vert
x\right\vert _{n}^{-1}x,\quad\left\vert \nabla u(x)\right\vert _{n}=\left\vert
\upsilon^{\prime}(\left\vert x\right\vert _{n})\right\vert \label{r-ngr}%
\end{equation}
for all $x\in\Delta$. Without any doubt, $u|_{\Omega}\in W^{1,\infty}%
(\Omega)$; moreover, $\nabla
u(x)\vec{n}=\upsilon^{\prime}(\left\vert x\right\vert
_{n})\left\vert x\right\vert _{n}^{-1}x\cdot(\left\vert x\right\vert
_{n}^{-1}x)=\upsilon^{\prime}(R_{1})$ for every $x\in\Sigma_{1}$ and
$\nabla u(x)\vec{n}=-\upsilon^{\prime}(R_{2})$ for $x\in\Sigma_{2}$.
Hence such a $u|_{\Omega}$ satisfies condition \cite[(7)]{LuGao:16}
if and only if
$\upsilon^{\prime}(R_{1})=\upsilon^{\prime}(R_{2})=0.$

Having in view the remark above, we discuss the result in \cite[Th.~1.1]%
{LuGao:16} for $W^{1,\infty}(\Omega)$ replaced by $C^{1}(\overline{\Omega}),$
more precisely the result in \cite{LuGao:16} concerning the local extrema of
$I$ defined in \cite[(1)]{LuGao:16} (quoted above) on the space
\begin{align*}
U:=  &  \{u:=\upsilon\circ\left\vert \cdot\right\vert _{n}\mid\upsilon\in
C^{1}[R_{2},R_{1}],~\upsilon^{\prime}(R_{1})=\upsilon^{\prime}(R_{2})=0\}\\
=  &  \{\upsilon\circ\left\vert \cdot\right\vert _{n}\mid\upsilon\in
C_{1,0}[R_{2},R_{1}]\}\subset C^{1}\left(  \overline{\Omega}\right)
\end{align*}
when $C^{1}\left(  \overline{\Omega}\right)  $ (and $U$) is endowed with the
norm
\begin{equation}
\left\Vert u\right\Vert :=\left\Vert u\right\Vert _{\infty}+\left\Vert \nabla
u\right\Vert _{\infty}; \label{r-gl13b}%
\end{equation}
moreover, in the sequel, $V:=C_{0}[R_{2},R_{1}]$ is endowed with the norm
$\left\Vert \cdot\right\Vert _{\infty}.$

Unlike \cite{LuGao:16}, let us set
\begin{equation}
F(r):=-\frac{1}{r^{n-1}}\int_{R_{2}}^{r}f(\rho)\rho^{n-1}d\rho=-\frac
{1}{r^{n-1}}\int_{R_{2}}^{r}\theta f,\quad r\in\lbrack R_{2},R_{1}],
\label{r-F}%
\end{equation}
where $\theta$ is defined in (\ref{r-gn-teta}).

\begin{remark}
\label{rem1}Notice that our $F(r)$ is $r$ times the one introduced in
\cite{LuGao:16}.
\end{remark}

From (\ref{r-F}) and the hypotheses on $f$, we have that
$F(R_{1})=F(R_{2})=0$ and $(\theta F)^{\prime}=-\theta f$ on
$[R_{2},R_{1}]$. Since
\[
(\theta F)^{\prime}(r)=0\iff(\theta f)(r)=0\iff f(r)=0\iff r=R_{3}%
\]
and $(\theta F)(R_{1})=(\theta F)(R_{2})=0$, it follows that $\theta
F>0$ or $\theta F<0$ on $(R_{2},R_{1})$, that is $F>0$ or $F<0$ on
$(R_{2},R_{1}).$ Moreover, from the definition of $F$ we get
\[
R_{2}^{n-1}\left\vert F(r)\right\vert \leq\left\vert r^{n-1}F(r)\right\vert
=\left\vert \int_{R_{2}}^{r}\theta f\right\vert \leq\int_{R_{2}}^{R_{1}%
}\left\vert \theta f\right\vert =\big\Vert\theta f\big\Vert_{L^{1}[R_{2}%
,R_{1}]}<R_{2}^{n-1}(2\lambda/3)^{3/2}%
\]
for every $r\in\lbrack R_{2},R_{1}]$, whence $\left\vert
F(r)\right\vert <(2\lambda/3)^{3/2}$ for $r\in\lbrack R_{2},R_{1}].$

Let $u\in U$, that is $u:=\upsilon\circ\left\vert \cdot\right\vert
_{n}$ with $\upsilon\in C_{1,0}[R_{2},R_{1}]$, and set
$v:=\upsilon^{\prime}$ $(\in
C_{0}[R_{2},R_{1}])$. We have that%
\begin{equation}
\int_{R_{2}}^{R_{1}}\theta f\upsilon=-\int_{R_{2}}^{R_{1}}(\theta F)^{\prime
}\upsilon=-(\theta F\upsilon)|_{R_{2}}^{R_{1}}+\int_{R_{2}}^{R_{1}}\theta
F\upsilon^{\prime}=\int_{R_{2}}^{R_{1}}\theta Fv. \label{r-gl14}%
\end{equation}

Using (\ref{r-int}) and (\ref{r-ngr}) we get%
\begin{align*}
I[u]  &  =\int_{\Omega}\left[  H(\left\vert \nabla u(x)\right\vert
)-f(\left\vert x\right\vert )u(x)\right]  dx=\int_{\Omega}\left[  H(\left\vert
\upsilon^{\prime}(\left\vert x\right\vert _{n})\right\vert )-f(\left\vert
x\right\vert )\upsilon(\left\vert x\right\vert )\right]  dx\\
&  =\gamma_{n}\int_{R_{2}}^{R_{1}}\theta(H\circ\left\vert v\right\vert
-Fv)=\gamma_{n}\int_{R_{2}}^{R_{1}}\theta(H\circ v-Fv),
\end{align*}
that is
\[
I[u]=\gamma_{n}K(v),
\]
where $K$ is defined in (\ref{r-K}) and $[a,b]:=[R_{2},R_{1}]$.
Therefore, Theorem \ref{t-gl} applies also in this situation.
Applying it we get that $I$ defined in \cite[(1)]{LuGao:16} has no
local minimizers and $\overline{u}\in C^{1}(\overline{\Omega})$ is a
local maximizer of $I|_{U}$ if and only if there exists
$u_{0}\in\mathbb{R}$ such that $\overline{u}(x)=u_{0}+\int
_{R_{2}}^{\left\vert x\right\vert _{n}}z_{2}(F(\rho))d\rho$ for
every $x\in\overline{\Omega}$, where $z_{2}(A)$ is the unique
solution of the equation $z\left(  \tfrac{1}{2}z^{2}-\lambda\right)
=A$ in the interval
$(-\sqrt{2\lambda/3},\sqrt{2\lambda/3})$ for $A\in(-(2\lambda/3)^{3/2}%
,(2\lambda/3)^{3/2}).$

For the present discussion we take the case in which $F>0$ on $(a,b)$. In this
case observe that $z_{2}(A)=A/E_{3}^{-1}(A^{2})$ for $A\in(0,(2\lambda
/3)^{3/2})$. This proves that the first and second assertions of \cite[Th.
1.1]{LuGao:16} are \emph{false}; in particular, $\overline{u}_{2}$ is not a
local minimizer of $I|_{U}$ (exactly as in the 1D case).

Moreover, from the discussion above, we can conclude that also the assertion
\textquotedblleft In comparison with the 1D case discussed by D. Gao and R.
Ogden, there exists huge difference in higher dimensions" from the abstract of
\cite{LuGao:16} is false.

\end{document}